\documentclass[11pt,twoside]{amsart}
\usepackage[latin1]{inputenc}
\usepackage[T1]{fontenc}
\usepackage{amsmath}
\usepackage{amsthm}
\usepackage{amssymb}
\usepackage[all]{xypic}

\newtheorem{thm}{Theorem}[section]

\newtheorem{lem}[thm]{Lemma}
\newtheorem{cor}[thm]{Corollary}

\numberwithin{equation}{section}

\theoremstyle{definition}
\newtheorem{definition}[thm]{Definition}
\newtheorem{remark}[thm]{Remark}

\newcommand\Tor{\mathcal T\!or}

\newcommand{\Pic}{{\rm Pic}}

\newcommand{\rk}{{\rm rk}}

\newcommand{\Hom}{{\rm Hom}}
\newcommand{\Spec}{{\rm Spec}}

\newcommand{\id}{{\rm id}}
\newcommand{\tr}{{\rm tr}}

\newcommand{\mono}{\hookrightarrow}

\newcommand{\Ext}{{\rm Ext}}

\newcommand{\sHom}{\cal{H}om}
\newcommand{\K}{{\rm K}}

\newcommand{\cal}{\mathcal}

\newcommand{\kf}{{\cal F}}
\newcommand{\kg}{{\cal G}}
\newcommand{\kh}{{\cal H}}

\newcommand{\kk}{{\cal K}}
\newcommand{\kl}{{\cal L}}

\newcommand{\km}{{\cal M}}
\newcommand{\ko}{{\cal O}}

\newcommand{\LL}{\mathbb{L}}

\newcommand{\CC}{\mathbb{C}}
\newcommand{\EE}{\mathbb{E}}

\newcommand{\HH}{\mathbb{H}}
\newcommand{\PP}{\mathbb{P}}

\renewcommand{\to}{\xymatrix@1@=15pt{\ar[r]&}}
\renewcommand{\rightarrow}{\xymatrix@1@=15pt{\ar[r]&}}
\renewcommand{\mapsto}{\xymatrix@1@=15pt{\ar@{|->}[r]&}}
\renewcommand{\twoheadrightarrow}{\xymatrix@1@=15pt{\ar@{->>}[r]&}}
\renewcommand{\hookrightarrow}{\xymatrix@1@=15pt{\ar@{^(->}[r]&}}
\newcommand{\congpf}{\xymatrix@1@=15pt{\ar[r]^-\sim&}}
\renewcommand{\cong}{\simeq}

\begin{document}

\title[Deformation-obstruction theory of complexes]{Deformation-obstruction
theory for complexes via Atiyah and Kodaira--Spencer classes}
\author[D.\ Huybrechts, R.\,P.\ Thomas]{Daniel Huybrechts and Richard Thomas}

\address{D.H.: Mathematisches Institut,
Universit{\"a}t Bonn, Beringstr.\ 1, 53115 Bonn, Germany}
\email{huybrech@math.uni-bonn.de}

\address{R.P.T.: Department of Mathematics, Imperial College, London SW7
2AZ, UK}\email{rpwt@ic.ac.uk}

\begin{abstract} \noindent
We give a universal approach to the deformation-obstruction theory
of objects of the derived category of coherent sheaves over a
smooth projective family. We recover and generalise the
obstruction class of Lowen and Lieblich, and prove that it is a
product of Atiyah and Kodaira--Spencer classes. This allows us to
obtain deformation-invariant virtual cycles on moduli spaces of
objects of the derived category on threefolds.
\vspace{-2mm}\end{abstract}

\maketitle

\section{Introduction}

The deformation theory of objects of the derived category of
coherent sheaves on a smooth projective family has been studied in
\cite{Lieb, Lowen}. These papers produce a class which is the
obstruction to deforming a complex sideways over an infinitesimal
deformation. The theory is not in its final form, however. The
obstruction class is expected to be a product of Atiyah and
Kodaira--Spencer classes, but this has only been proved in some
situations (such as \cite{HMS}). As a consequence the results are
not strong enough to obtain virtual cycles \cite{BF, LT}, as
discussed in \cite[Sect.\ 2]{PT}, for instance.

This paper generalises \cite{Lieb, Lowen} by working universally,
using Fourier--Mukai kernels as in \cite{HMS, HT}. We put the
discussion of \cite[App.\ C]{HMS}, which deals with deformations
over $\Spec(k[t]/t^n)$, in a more general context. This allows us
to identify the obstruction class of \cite{Lieb, Lowen} with the
product of certain Atiyah and Kodaira--Spencer classes. We are
then able to obtain virtual cycles on moduli spaces of simple
complexes on threefolds, completing the foundational work required
in \cite[Sect.\ 2]{PT}.

While the classical Atiyah and Kodaira--Spencer classes were enough for the
case treated in \cite{HMS}, Illusie's versions of these classes, defined
in \cite{Ill} using the cotangent complex, are necessary to go beyond the
simplest cases. This was done for complexes of modules in \cite{Ill} and
for a single sheaf in \cite{BuFl1}.
However, as with the theory of virtual cycles \cite{BF, LT},
only a small part of the cotangent complex affects the
computations. For this reason we work with truncated versions of the cotangent
complex, Atiyah class and Kodaira--Spencer class. These can be written
down simply and explicitly, without reference to Illusie's versions, after
picking an embedding into a smooth ambient space. This keeps our presentation
elementary, avoiding dg and simplicial resolutions.


Suppose there exists an embedding of a noetherian separated
scheme $X$ in a smooth ambient space $A$ with ideal sheaf
$J\subset \ko_A$. Then the truncation $\tau^{\geq-1}L^\bullet_X$ of Illusie's
cotangent complex is quasi-isomorphic to the complex
$$
\LL_X:=(J/J^2\xymatrix{\ar[r]&}\Omega_A|_X)
$$
concentrated in degrees $-1$ and $0$ (see Section
\ref{sect:ATKS}).
 Then
for any perfect complex $E$ on $X$ we introduce the truncated Atiyah
class
$$
A(E)\in\Ext^1_X(E,E\otimes \LL_X).
$$
The class $A(E)$ lifts the classical Atiyah class in
$\Ext^1_X(E, E\otimes\Omega_X)$ via the obvious map
$\LL_X\to\Omega_X$. Moreover, it can be shown that
$A(E)$ is the truncation of Illusie's Atiyah class in
$\Ext^1_X(E,E\otimes L^\bullet_X)$, though we will not
need this. \smallskip

Suppose that $i\colon X_0\
\hookrightarrow\,X$ is a square zero thickening of schemes, i.e.
$i$ is a closed embedding defined by an ideal
sheaf $I$ on $X$ with $I^2=0$. We inherit an embedding $X_0\subset A$
from $X\subset A$, yielding
a truncated cotangent complex $\LL_{X_0}$.
Then one can also define a
truncated Kodaira--Spencer class
$$\kappa(X_0/X)\in\Ext^1_{X_0}(\LL_{X_0},I)$$
of the embedding $X_0\subset X$. Given a perfect complex
$E_0$ over $X_0$, there is, as above, a truncated Atiyah class $A(E_0)\in\Ext^1_{X_0}
(E_0,E_0\otimes\LL_{X_0})$.
The main result of this paper is that the product of these classes
is the obstruction to deforming the complex $E_0$ over $X$. 
\medskip

\noindent{\bf Theorem} {\it Let $E_0$ be a perfect complex on a
separated noetherian scheme $X_0$ and let $i\colon X_0\
\hookrightarrow\,X$ be a closed embedding defined by an ideal $I$
of square zero. Assume that $X$ can be embedded into a smooth
ambient space $A$.  Then there exists a perfect complex
$E$ on $X$ such that the derived pull-back $i^*E$ is
quasi-isomorphic to $E_0$ if and only if
$$0=({\rm id}_{E_0}\otimes\kappa(X_0/X))\circ A(E_0)\ \in\
\Ext^2_{X_0}(E_0,E_0\otimes I).$$}

The existence of the obstruction class is already proved in \cite{Lieb, Lowen}
(in the slightly less general product situation of $X_0=X\times_{\Spec(A)}\Spec(A_0)$
). But the explicit expression as
the product of the truncated Atiyah class of
$E_0$ and the truncated Kodaira--Spencer class of the thickening $X_0\subset
X$ is crucial for some applications \cite{HMS, PT}.

Although there is also a version of the theorem phrased
in terms of the full cotangent complex $L^\bullet_X$ and Illusie's corresponding
Atiyah and Kodaira--Spencer classes, only the degree $-1$ and $0$ parts of
$L^\bullet_X$ contribute to the obstruction because the lower degree
terms yield classes in $\Ext^{\ge3}$.

\medskip

In this paper we emphasize the view that all classes (Atiyah,
Kodaira--Spencer, obstruction) are most naturally defined
universally. Suppose we have a class $a(E)\in \Ext^k_X(E\otimes F_1,E\otimes
F_2)$ which is natural in $E$ for fixed complexes $F_1$ and $F_2$ on $X$.
Thinking of it as a natural transformation of functors
$$
(\ \cdot\ )\otimes F_1\xymatrix{\ar[r]&}(\ \cdot\ )\otimes F_2[k]
$$
then by the philosophy of Fourier--Mukai transforms one expects it to be
represented by a morphism of complexes
$$
a\colon i_{\Delta_X*}(F_1)\xymatrix{\ar[r]&} i_{\Delta_X*}(F_2)[k]
\quad\rm{on}\ X\times X.
$$
That is $a(E)$ should be obtained by tensoring $a$ with $\pi_1^*E$ and then
taking the direct image $\pi_{2*}$. (Here $\pi_i$ is projection onto the
$i$th factor of $X\times X$ and $i_{\Delta_X}$ is the diagonal inclusion.)

So in this paper we define \emph{universal} truncated Atiyah and Kodaira--Spencer
classes, and their product defines a \emph{universal} obstruction class
$$
\varpi\colon \ko_{\Delta_{X_0}}\xymatrix{\ar[r]&}
i_{\Delta_{X_0}*}(I)[2].
$$
Applying this to any complex $E_0$ yields the obstruction class
$\varpi(E_0)$ described by the main theorem. This avoids the ad-hoc constructions
of Lieblich. Note that universal Atiyah classes
were studied via the full cotangent complex in \cite{BuFl2}, but for the
truncated version our approach is much more elementary.

A topic of current interest \cite{Joy, KS, PT, Th, Toda} is
defining invariants of threefolds $X$ (and particularly
Calabi--Yau threefolds) using virtual cycles on moduli spaces of
objects of $D^b(X)$ satisfying some natural conditions (such as
stability or simplicity). This paper fills the gap mentioned in
\cite[Sect.\ 2]{PT} to give such virtual cycles.
\medskip

{\bf Notation.}
%
All our schemes will be separated and noetherian over a fixed field $k$.
One could work relative to a fixed base, but then various flatness issues have to be addressed.

In the sheaf-theoretic computations of Section \ref{sect:ATKS}
most of our functors are the usual ones; in Sections
\ref{sect:Lieblich} and \ref{sect:appl} they are all derived.
Thus, for instance, Hom denotes the (quasi-isomorphism class of)
the complex RHom, with cohomologies $\kh^i(\Hom)=\Ext^i$.

We use $I_{Y\subset Z}$ to denote the ideal sheaf
of a subscheme $Y$ of another
scheme $Z$, with a few exceptions. We reserve $I$ to denote the square zero
ideal of $i\colon X_0\subset X$, and $J$ denotes the ideal of an ambient
embedding of $X$ in some smooth $A$. The ideal of the induced embedding $X_0\subset
A$ is denoted by $J_0$. And for any separated
scheme $Y$, we use $I_{\Delta_Y}$ to denote the ideal of the diagonal $\Delta_Y$,
the image of the diagonal embedding $i_{\Delta_Y}\colon Y\to Y\times Y$.

We sometimes suppress pushforwards by embeddings, to keep the notation
to a minimum. Thus $\ko_{\Delta_Y}$ might denote $i_{\Delta_Y*}\ko_Y$.
\medskip

{\bf Acknowledgements.} We would like to thank Max
Lieblich, Emanuele Macr\`i and Paolo Stellari for valuable comments on a
first version of this paper, Fabian Langholf and Pierrick Bousseau for pointing out a problem with flatness in the second version, and Barbara Fantechi,
Ian Grojnowksi and Yinan Song for useful remarks. Financial
support for the first author by Imperial College, EPSRC and
SFB/TR45 of the DFG are gratefully acknowledged. The second author
thanks the Royal Society and Leverhulme Trust for support.

\section{Truncated Atiyah and Kodaira--Spencer
classes}\label{sect:ATKS}

\subsection{The truncated cotangent
complex}\label{subsect:trunccot} Fix a closed embedding $X\subset
A$ of a scheme $X$ into a smooth scheme $A$ and let
$J\subset\ko_A$ be its ideal sheaf.

The composition of the differential $d_A\colon J\subset
\ko_A\to\Omega_A$ with the
projection $\Omega_A\to \Omega_A|_X$ is
an $\ko_A$-module homomorphism
. It factors through $J/J^2$ to become
an $\ko_X$-linear homomorphism
$$d_{X/A}\colon J/J^2\xymatrix{\ar[r]&}\Omega_A|_X$$
with cokernel $\Omega_X$.

\begin{definition}
The \emph{truncated cotangent complex} of $X$ (in  $A$) is the
length-two complex
$$
\spreaddiagramcolumns{1pc}
\LL_X:=\xymatrix{\Big(J/J^2\ar[r]^{d_{X/A}} &
\Omega_A|_X\Big)}
$$ concentrated in degrees $-1$ and $0$. By construction, $\LL_X$ comes with
a natural homomorphism $\LL_X\to\kh^0(\LL_X)\cong\Omega_X$.
\end{definition}

As the notation suggests, the complex $\LL_X$ is independent (up
to quasi-isomorphism) of the embedding $X\subset A$. This is proved by
showing that $\LL_X$ is quasi-isomorphic to the truncation
$\tau^{\geq-1}(L^\bullet_X)$ of Illusie's cotangent complex, or directly
as follows. Suppose that $X$ embeds into two smooth
schemes $A_1,\,A_2$ with ideals $J_1,\,J_2$, giving the diagonal embedding
in $A_1\times A_2$ with ideal $J_{12}$. Then the short exact sequence of
2-term complexes
\begin{equation} \label{indep}
\spreaddiagramrows{-0.6pc}
\xymatrix{
\big(J_1/J_1^2 \ar[r]\ar[d] & \Omega_{A_1}|_X\big) \ar[d] \\
\big(J_{12}/J_{12}^2 \ar[r]\ar[d] & \Omega_{A_1}|_X\oplus\Omega_{A_2}|_X
\big) \ar[d] \\
\big(\Omega_{A_2}|_X \ar@{=}[r] & \Omega_{A_2}|_X\big), }
\end{equation}
induced by the composition $X\subset A_1\times A_2\to A_1$,
shows that $\LL_X$ defined via $A_1$ is quasi-isomorphic to $\LL_X$ defined
via $A_1\times A_2$.

\subsection{The truncated Atiyah class}

The following simple Lemma is really the key to this paper. It does not require
the assumption that $A$ is smooth.

\begin{lem} \label{lem:k}
For any $X\subset A$ with ideal sheaf $J\subset\ko_A$, the kernel of the
natural surjection $I_{\Delta_A}|_{X\times X}\to I_{\Delta_X}$ is isomorphic
to
$$
\Tor_1^{X\times A}(\ko_{\Delta_X},\ko_{X\times X})\cong i_{\Delta_X*}(J/J^2).
$$
\end{lem}

\begin{proof}
We first claim that the restriction of the sheaf $I_{\Delta_A}$ to $X\times
A$ is isomorphic to $I_{\Delta_X\subset X\times A}$.
Tensoring the exact sequence
$$
0\xymatrix{\ar[r]&} J\boxtimes\ko_A\xymatrix{\ar[r]&}\ko_{A\times
A}\xymatrix{\ar[r]&}\ko_{X\times A}\xymatrix{\ar[r]&}0
$$
by $\ko_{\Delta_A}$ gives the exact sequence
$$
i_{\Delta_A*}J\xymatrix{\ar[r]&}\ko_{\Delta_A}\xymatrix{\ar[r]&}\ko_{\Delta_X}\xymatrix{\ar[r]&}0.
$$
Since the leftmost arrow is an injection, $\Tor_1^{A\times A}(\ko_{X\times
A},\ko_{\Delta_A})=0$. Therefore restricting $0\to I_{\Delta_A}\to\ko_{A\times
A}\to\ko_{\Delta_A}\to0$ to $X\times A$ is exact:
\begin{equation} \label{eqn:diagA}
0\xymatrix{\ar[r]&} I_{\Delta_A}|_{X\times
A}\xymatrix{\ar[r]&}\ko_{X\times
A}\xymatrix{\ar[r]&}\ko_{\Delta_X}\xymatrix{\ar[r]&}0.
\end{equation}
Thus indeed $I_{\Delta_A}|_{X\times A}\cong I_{\Delta_X\subset X\times
A}$.

Restricting (\ref{eqn:diagA}) to $X\times X$ yields the exact
sequence
\begin{equation} \label{exlater}
\xymatrix{0\ar[r]&\kk\ar[r]&I_{\Delta_A}|_{X\times X}\ar[r]& \ko_{X\times
X}\ar[r]&\ko_{\Delta_X}^{\phantom{I^I}}\ar[r]&0,}
\end{equation}
where $\kk$ is the kernel we want to compute. Therefore
$\kk\cong\Tor_1^{X\times A}(\ko_{\Delta_X},\ko_{X\times X})$.
Computing this by tensoring the exact sequence
\begin{equation} \label{late}
\xymatrix{0\ar[r]&\ko_X\boxtimes J\ar[r]&\ko_{X\times
A}\ar[r]&\ko_{X\times X}\ar[r]&0}
\end{equation}
with $\ko_{\Delta_X}$ gives
\begin{equation} \label{kk}
\kk\cong\Tor_1^{X\times A}(\ko_{\Delta_X},\ko_{X\times
X})\cong(\ko_X\boxtimes J)\otimes\ko_{\Delta_X}\cong
i_{\Delta_X*}(J/J^2). \vspace{-19pt}
\end{equation}
\end{proof}

Thus we have a quasi-isomorphism $\ko_{\Delta_X}\cong
\big(i_{\Delta_X*}(J/J^2)\to I_{\Delta_A}|_{X\times X}\to\ko_{X\times X}\big)$,
with the complex in degrees $-2,\,-1$ and 0.

\begin{definition} \label{def:At} The \emph{truncated universal Atiyah class}
$$\alpha_X\in\Ext^1_{X\times X}(\ko_{\Delta_X},i_{\Delta_X*}(\LL_X))$$ is
given by the map of complexes
\begin{equation} \label{Atcx}
\spreaddiagramrows{-.8pc}
\xymatrix{
\ko_{\Delta_X}\ \cong & \Big(i_{\Delta_X*}(J/J^2) \ar[r]\ar@{=}[d] & I_{\Delta_A}|_{X\times
X} \ar[r]\ar[d] & \ko_{X\times X}\Big) \\
i_{\Delta_X*}\LL_X[1]\ \cong & \Big(i_{\Delta_X*}(J/J^2) \ar[r] & (I_{\Delta_A}/I_{\Delta_A}^2)
\big|^{}_{X\times X}\Big), \hspace{-1cm}}
\end{equation}
where the second vertical arrow is the obvious projection. For this to make
sense we must first check that the diagram \eqref{Atcx}
commutes. Define $\varphi\colon I_{\Delta_A}|_{X\times A}\to\ko_{X\times
X}$ by
$$
\xymatrix{
&0\ar[r]&I_{\Delta_A}|_{X\times A}\ar@{->>}[d]\ar[dr]_{\varphi}\ar[r]&\ko_{X\times
A}\ar@{->>}[d]\ar[r]&\ko_{\Delta_X}\ar@{=}[d]\ar[r]&0\\
0\ar[r]&i_{\Delta_X*}(J/J^2)\ar[r]&I_{\Delta_A}|_{X\times X}\ar[r]&\ko_{X\times
X}\ar[r]&\ko_{\Delta_X}\ar[r]&0.}
$$
In the bottom row we are suppressing the pushforward map from $X\times X$
to $X\times A$. The kernel of $\varphi$ is $(I_{\Delta_X\subset X\times A})\cap(\ko_X\boxtimes
J)=\ko_X\boxtimes J$, which yields
$$\xymatrix{\ko_X\boxtimes J\ar@{}[rrd]|\circlearrowleft\ar@{->>}[d]_{\psi}\ar[rr]&&
I_{\Delta_A}|_{X\times A}\ar@{->>}[d] \\
i_{\Delta_X*}(J/J^2)\ar[rr]&&I_{\Delta_A}|_{X\times X}.}$$

Since the inclusion of $\ko_X\boxtimes J$ inside $I_{\Delta_A}|_{X\times
A}\subset\ko_{X\times A}$ is part of the exact sequence \eqref{late}, the
description \eqref{kk} of $i_{\Delta_X*}(J/J^2)$
as $(\ko_X\boxtimes J)\otimes\ko_{\Delta_X}$ shows that $\psi$ is the natural
projection $\ko_X\boxtimes J\to i_{\Delta_X*}(\ko_X\boxtimes J)|_{\Delta_X}\cong
i_{\Delta_X*}(J/J^2)$. So it suffices to show that
$$\xymatrix{
\ko_X\boxtimes J\ar[d]\ar[rr]&&I_{\Delta_A}|_{X\times A}\ar[d] \\
i_{\Delta_X*}(J/J^2)\ar[rr]^{i_{\Delta_X*}(d_{X/A})} &&
(I_{\Delta_A}/I_{\Delta_A}^2)\big|^{}_{X\times X}\hspace{-1cm}}
$$
commutes. Given $f\in J$, $d_{X/A}(f)$ is the image of $1\otimes f-f\otimes
1$ in $(I_{\Delta_A}/I_{\Delta_A}^2)|^{}_{X\times X}$. But
$(f\otimes 1)|_{X\times A}=0$,
so this is the image of $1\otimes f\in \ko_X\boxtimes
J\subset I_{\Delta_A}|_{X\times A}$ in $(I_{\Delta_A}/I_{\Delta_A}^2)
|^{}_{X\times X}$.
\end{definition}

Using methods like those of Section  \ref{subsect:trunccot} one
can also show that the truncated Atiyah class is independent of
the embedding $X\subset A$.

Composing \eqref{Atcx} with the direct image of the natural map
$\LL_X\to\kh^0(\LL_X)$ gives
$$
\spreaddiagramrows{-.8pc}
\xymatrix{
\Big(I_{\Delta_X} \ar[r]\ar[d] & \ko_{X\times X}\Big) \\
I_{\Delta_X}/I_{\Delta_X}^2.\!\!}
$$
(Here we have used the isomorphism $(I_{\Delta_A}|_{X\times X})\big/i_{\Delta_X*}(J/J^2)
\cong I_{\Delta_X}$ of Lemma \ref{lem:k}.) Dividing the upper row by the
acyclic complex $I_{\Delta_X}^2\to I_{\Delta_X}^2$ we get the class in
$\Ext^1_{X\times X}(\ko_{\Delta_X},i_{\Delta_X*}\Omega_X)$ of the extension
\begin{equation} \label{class}
0\xymatrix{\ar[r]&}
i_{\Delta_X*}\Omega_X\xymatrix{\ar[r]&}\ko_{X\times
X}/I_{\Delta_X}^2\xymatrix{\ar[r]&}\ko_{\Delta_X}\xymatrix{\ar[r]&}0.
\end{equation}
Therefore $\alpha_X$ projects to the classical universal Atiyah class.

\begin{remark}
Another way of phrasing the construction is as follows. The natural map
\begin{equation} \label{projI}
I_{\Delta_X}\xymatrix{\ar[r]&} i_{\Delta_X*}\Omega_X
\end{equation}
(given by the projection $I_{\Delta_X}\to I_{\Delta_X}/I^2_{\Delta_X}$) has
a natural lift
\begin{equation} \label{liftL}
I_{\Delta_X}\xymatrix{\ar[r]&} i_{\Delta_X*}\LL_X
\end{equation}
given by the map of complexes
$$
\spreaddiagramrows{-.8pc}
\xymatrix{
\Big(i_{\Delta_X*}(J/J^2) \ar[r]\ar@{=}[d] & I_{\Delta_A}|_{X\times X}\big)
\ar[d]  & \cong\ I_{\Delta_X} \\
\Big(i_{\Delta_X*}(J/J^2) \ar[r] &
(I_{\Delta_A}/I_{\Delta_A}^2)\big|^{}_{X\times X}\Big) &\cong\
i_{\Delta_X*}\LL_X\,.}
$$
Consider the boundary map of exact sequence $0\to I_{\Delta_X}\to\ko_{X\times
X}\to\ko_{\Delta_X}\to0$:
$$
\ko_{\Delta_X}\xymatrix{\ar[r]&} I_{\Delta_X}[1].
$$
The classical Atiyah class is its composition with \eqref{projI},
while the truncated Atiyah class is its composition with
\eqref{liftL}.
\end{remark}

\begin{remark}
After pushing forward $\alpha_X$ by the inclusion
$\iota\colon X\times X\ \hookrightarrow\,X\times A$, we claim that it can
be described by the obvious diagram
\begin{equation} \label{jjdg}
\spreaddiagramrows{-1pc}
\spreaddiagramcolumns{-1pc}
\xymatrix{
& J/J^2 \ar@{=}[r]\ar[d] & J/J^2 \ar[d] \\
0 \ar[r] & \Omega_A|_X \ar[r]\ar[d]
& \ko_{X\times A}/I_{\Delta_X\subset X\times A}^2 \ar[r]\ar[d] & \ko_{\Delta_X}
\ar[r]\ar@{=}[d] & 0 \\
0 \ar[r] & \Omega_X \ar[r]\ar[d] &
\ko_{X\times X}/I_{\Delta_X}^2 \ar[r]\ar[d] & \ko_{\Delta_X} \ar[r]
& 0.\! \\ & 0 & 0}
\end{equation}
Here we have suppressed some pushforward maps for appearance's sake.
The middle row is the restriction to $X\times A$ of $0\to I_{\Delta_A}/I_{\Delta_A}^2
\to\ko_{A\times A}/I_{\Delta_A}^2\to \ko_{\Delta_A}\to0$, using the vanishing
of $\Tor_1^{A\times A}(\ko_{X\times
A},\ko_{\Delta_A})$ (as shown in the proof of Lemma \ref{lem:k}).

The bottom row is (the pushforward of) the extension \eqref{class} defining
the classical Atiyah class. We view the top two rows as a horizontal short
exact sequence of two-term vertical complexes, giving
an extension class in $\Ext^1_{X\times A}(\ko_{\Delta_X},\iota_*i_{\Delta_X*}
(\LL_X))$ projecting to the classical Atiyah class. Write this extension
class as
\begin{equation} \label{ialpha}
\spreaddiagramrows{-1pc}
\spreaddiagramcolumns{-1pc}
\xymatrix{
&& \!\!\big(\iota_*i_{\Delta_X*}\Omega_A|_X \ar[r]\ar@{=}[d] & \ko_{X\times
A}/I_{\Delta_X\subset X\times A}^2\big) & \cong\ \ko_{\Delta_X} \\
\iota_*i_{\Delta_X*}\LL_X[1]\ \cong & \big(\iota_*i_{\Delta_X*}J/J^2 \ar[r]
& \iota_*i_{\Delta_X*}\Omega_A|_X\big)\,.\!\!\!\!}
\end{equation}
This is quasi-isomorphic to
$$
\spreaddiagramrows{-1pc}
\spreaddiagramcolumns{-1pc}
\xymatrix{
&& \big(I_{\Delta_X\subset X\times A} \ar[r]\ar[d] & \ko_{X\times
A}\big) & \cong\ \ko_{\Delta_X} \\
\iota_*i_{\Delta_X*}\LL_X[1]\ \cong & \big(\iota_*i_{\Delta_X*}J/J^2 \ar[r]
& \iota_*i_{\Delta_X*}\Omega_A|_X\big)\,,}
$$
as can be seen by dividing the two sheaves on the top row by $I_{\Delta_X\subset
X\times A}^2$. Using \eqref{eqn:diagA}, the vertical arrow is the projection
$I_{\Delta_A}|_{X\times A}\to(I_{\Delta_A}/I_{\Delta_A}^2)|_{X\times A}$.

Since any sheaf $F$ on $X\times A$ has a natural map to $F\otimes\iota_*\ko_{X\times
X}=\iota_*(F|_{X\times X})$, this maps term by term  and quasi-isomorphically
to the map of complexes
$$
\spreaddiagramrows{-1pc}
\spreaddiagramcolumns{-1pc}
\xymatrix{
& \big(\iota_*i_{\Delta_X*}J/J^2 \ar[r]\ar@{=}[d] & \iota_*(I_{\Delta_A}
|_{X\times X}) \ar[r]\ar[d] & \iota_*\ko_{X\times X}\big) & \cong\
\ko_{\Delta_X} \\
\iota_*i_{\Delta_X*}\LL_X[1]\ \cong \vspace{-1cm} & \big(\iota_*i_{\Delta_X*}J/J^2
\ar[r] & \iota_*i_{\Delta_X*}\Omega_A|_X\big)\,.}
$$
But this is precisely the pushforward of our
original definition \eqref{Atcx} of $\alpha_X$.
\end{remark}

\begin{definition} \label{rem:truncAt}
Thinking of \eqref{Atcx} as a map of Fourier--Mukai kernels, we apply it
to a perfect complex $E$ on $X$ to yield the \emph{truncated Atiyah
class} of $E$,
\begin{equation} \label{eq:truncAt}
A(E):=\alpha_X(E)\in\Ext^1_X(E,E\otimes\LL_X),
\end{equation}
mapping under $\LL_X\to\Omega_X$ to the classical Atiyah class in $\Ext^1_X(E,E\otimes\Omega_X)$.
\end{definition}

\subsection{Square zero extensions and Kodaira--Spencer classes} Now fix $i\colon X_0\ \hookrightarrow\,X$, a closed embedding defined by an ideal sheaf $I$ of square zero: $I^2=0$.
This allows us to consider $I$ to be an $\ko_{X_0}$-module. From $X\subset A$ we inherit a closed embedding $X_0\subset A$ with ideal $J_0$. The natural map $J_0/J_0^2\to I/I^2=I$ yields a morphism of complexes
\begin{equation} \label{dg:KS} \spreaddiagramrows{-.8pc} \xymatrix{
J_0/J_0^2 \ar[r]\ar[d] & \Omega_A|_{X_0} \\ I\,,\!} \end{equation}
where the upper row is the (shifted) truncated cotangent complex $\LL_{X_0}$ of $X_0$.

\begin{definition}\label{def:truncunivKS}
The \emph{truncated Kodaira--Spencer class} of the square zero extension $i\colon X_0\ \hookrightarrow\,X$ is the extension class $$\kappa({X_0/X})\in\Ext^1_{X_0}(\LL_{X_0},I)$$ defined by (\ref{dg:KS}).
\end{definition}

\begin{remark}
The map of complexes \eqref{dg:KS} factors through the map of complexes \begin{equation} \label{dg:KS1} \spreaddiagramrows{-.8pc} \xymatrix{ I \ar[r]\ar[d] & \Omega_X|_{X_0} \\ I\,.\!} \end{equation} When $I\to\Omega_X|_{X_0}$ is an injection (i.e. when it is quasi-isomorphic as a complex to $\Omega_{X_0}$), \eqref{dg:KS1} is the \emph{classical} Kodaira--Spencer class in $\Ext^1_{X_0}(\Omega_{X_0},I)$ and the truncated class is its pullback via the map $\LL_{X_0}\to\Omega_{X_0}$.
\end{remark}

By Lemma \ref{lem:k} applied to both $(X_0\subset A,J_0)$ and $(X_0\subset X,I)$ in place of $(X\subset A,J)$ we get the two horizontal short exact sequences \begin{equation}\label{eqn:sesXX0}
\spreaddiagramrows{-.5pc}
\xymatrix{
0\ar[r]& i_{\Delta_{X_0}*}(J_0/J_0^2) \ar[r]\ar[d] & I_{\Delta_A}|_{X_0\times X_0} \ar[r]\ar[d] & I_{\Delta_{X_0}} \ar[r]\ar@{=}[d] &0 \\ 0\ar[r]& i_{\Delta_{X_0}*}(I)\ar[r]& I_{\Delta_X}|_{X_0\times X_0}\ar[r]& I_{\Delta_{X_0}}\ar[r]&0\,.\!\!} \end{equation} The first vertical arrow is induced by the other two. The top left hand corner, composed with the surjection $I_{\Delta_A}|_{X_0\times X_0}\to (I_{\Delta_A}/I_{\Delta_A}^2)\big|^{}_{X_0\times X_0}$, gives \begin{equation}\label{i*KS} \spreaddiagramrows{-.8pc} \xymatrix{
i_{\Delta_{X_0}*}(J_0/J_0^2) \ar[r]\ar[d] & (I_{\Delta_A}/I_{\Delta_A}^2) \big|^{}_{X_0\times X_0} \\ i_{\Delta_{X_0}*}(I).\!}
\end{equation}
This is the pushforward of the diagram \eqref{dg:KS} by $i_{\Delta_{X_0}*}$.

\subsection{Product of the two classes}

It follows that the product of the universal truncated Atiyah
class \eqref{Atcx} of $X_0$ and the pushforward by
$i_{\Delta_{X_0}*}$ of the universal truncated Kodaira--Spencer
class of $X_0\subset X$, which is \eqref{i*KS}, can be represented
by the morphism of complexes
\begin{equation} \label{product}
\spreaddiagramrows{-.5pc}
\xymatrix{
\!\!\Big(i_{\Delta_{X_0}*}(J_0/J_0^2) \ar[r]\ar[d] &
I_{\Delta_A}|_{X_0\times X_0} \ar[r] & \ko_{X_0\times X_0}\Big) \\
i_{\Delta_{X_0}*}(I)\,.\!\!}
\end{equation}
Consider \eqref{eqn:sesXX0} to be a (vertical) quasi-isomorphism between
two 2-term (horizontal) complexes that are themselves quasi-isomorphic to
$I_{\Delta_{X_0}}$. By construction the above vertical map factors through
this quasi-isomorphism. Therefore \eqref{product} can be written
\begin{equation} \label{product2}
\spreaddiagramrows{-.5pc}
\xymatrix{
\!\!\Big(i_{\Delta_{X_0}*}(I) \ar[r]\ar[d] & I_{\Delta_X}|_{X_0\times
X_0} \ar[r] & \ko_{X_0\times X_0}\Big) \\
i_{\Delta_{X_0}*}(I)\,,\!\!}
\end{equation}
with vertical arrow the identity. That is, $i_{\Delta_{X_0}*} (\kappa(X_0/X))\circ\alpha_{X_0}
\in\Ext^2_{X_0\times X_0}(\ko_{\Delta_{X_0}},
i_{\Delta_{X_0}*}(I))$ is the extension class of the exact
sequence
\begin{equation} \label{product3}
0\xymatrix{\ar[r]&} i_{\Delta_{X_0}*}(I)\xymatrix{\ar[r]&}
I_{\Delta_X}|_{X_0\times X_0}\xymatrix{\ar[r]&}\ko_{X_0\times
X_0}\xymatrix{\ar[r]&}\ko_{\Delta_{X_0}}\xymatrix{\ar[r]&}0.
\end{equation}

\subsection{The universal obstruction class}

Let $h=\id\times i$ denote the natural closed embedding
$X_0\times X_0\ \hookrightarrow\,X_0\times X$ and define
$$H:=h^*h_*\ko_{\Delta_{X_0}},$$ where here $h^*$ is the \emph{derived}
pull-back. Then $H$ is a complex concentrated in degree $\leq0$.
Now $\Tor_0^{X_0\times X}(\ko_{\Delta_{X_0}},\ko_{X_0\times
X_0})\cong\ko_{\Delta_{X_0}}$, and $\Tor_1^{X_0\times
X}(\ko_{\Delta_{X_0}},\ko_{X_0\times X_0})\cong
i_{\Delta_{X_0}*}(I)$ by Lemma \ref{lem:k} applied to $(X_0\subset
X,I)$. Therefore
\begin{equation}\label{eqn:coh0-1}
\kh^0(H)\cong\ko_{\Delta_{X_0}}~~\text{and}~~\kh^{-1}(H)\cong i_{\Delta_{X_0}*}(I).
\end{equation}
This leads to the following definition (whose name will be justified in Section
\ref{sect:Lieblich}).

\begin{definition} \label{uoc} The \emph{universal obstruction class}
$$
\varpi(X_0/X)\in\Ext^2_{X_0\times X_0}(\ko_{\Delta_{X_0}},i_{\Delta_{X_0}*}(I))
$$
of $X_0\subset X$ is the extension class of the exact triangle
\begin{equation}\label{tauH}
\kh^{-1}(H)[1]\xymatrix{\ar[r]&}\tau^{\ge-1}(H)\xymatrix{\ar[r]&}\kh^0(H).
\end{equation}
\end{definition}

\begin{thm} \label{product4}
$\varpi(X_0/X)$ equals the product $i_{\Delta_{X_0}*}(\kappa(X_0/X))
\circ\alpha_{X_0}$.
\end{thm}

\begin{proof}
Consider the exact sequence
\begin{equation} \label{resol}
0\xymatrix{\ar[r]&} I_{\Delta_{X_0}\subset X_0\times
X}\xymatrix{\ar[r]&}\ko_{X_0\times X}\xymatrix{\ar[r]&}
h_*\ko_{\Delta_{X_0}}\xymatrix{\ar[r]&}0
\end{equation}
as a resolution of $h_*\ko_{\Delta_{X_0}}$. Since $\ko_{X_0\times
X}$ is flat, the result of applying the \emph{un}derived functor $h^*$ to
the resolution is $\tau^{\ge-1}(H)$. So we
restrict to $X_0\times X_0$, giving
\begin{equation} \label{t-H}
\spreaddiagramcolumns{-.5pc}
\spreaddiagramrows{-.8pc}
\xymatrix{
0 \ar[r] & \Tor_1^{X_0\times X}(\ko_{\Delta_{X_0}},\ko_{X_0\times X_0}) \ar[r]
\ar@{=}[d] & I_{\Delta_X}|_{X_0\times X_0} \ar[r] & \ko_{X_0\times X_0} \ar[r]
& \ko_{\Delta_{X_0}}
\ar[r] & 0, \\ & i_{\Delta_{X_0}*}(I)}
\end{equation}
by Lemma \ref{lem:k}. The two terms in the middle represent the
complex $\tau^{\ge-1}(H)$, and the exact sequence represents the
triangle \eqref{tauH}.

But \eqref{t-H} is precisely the exact sequence \eqref{product3}: both
are by construction the exact sequence \eqref{exlater} of Lemma \ref{lem:k},
applied to $(X_0\subset X,I)$ in place of $(X\subset A,J)$.
\end{proof}

\begin{remark} i) Recall that when $I$ injects into $\Omega_X|_{X_0}$, the class
$\kappa(X_0/X)$
is the pullback of the classical Kodaira--Spencer class \eqref{dg:KS1}. In
this case one can show that the above product is just the product of the
\emph{classical} Atiyah and Kodaira--Spencer classes.

This is the case, for instance, when $X$ is flat over
$\Spec(k[t]/t^{n+1})$ and $X_0$ is its base-change to $\Spec(k[t]/t^n)$:
the situation studied in \cite{HMS}. This explains why it was not necessary
to work with (truncated) cotangent complexes in that paper.

ii) In \cite{HMS} the base field is $k=\CC$ and $X_0$
is deformed into a non-algebraic direction. In particular the
situation studied there is one of the rare examples where $X$
cannot be embedded into a smooth ambient scheme $A$ over
$\Spec(\CC)$. However, an embedding into a smooth formal scheme
can still be found. Indeed, in \cite{HMS} $X$ is the $n$th order
neighbourhood of a smooth holomorphic family over a
one-dimensional disk and we can let $A$ be the formal neighbourhood of the
special fibre inside
this holomorphic family. The arguments in this section work
equally well in this setting.
\end{remark}

\subsection{Relative version.}\label{sec:relKS}
In the application in Section \ref{sect:appl} we shall need a relative version of the Kodaira--Spencer class.
Assume $X$ comes with a morphism $X\to B$ and the embedding $X\subset A$ is part of a commutative diagram
$$
\xymatrix{X_0\,\ar@{^(->}[r]&X\,\ar[dr]\ar@{^(->}[r]&A_B\,\ar@{^(->}[r]
\ar[d]&A\ar[d] \\
&&B\ \ar@{^(->}[r]&\widetilde B,\!}
$$
with $\widetilde B$ and $A\to \widetilde B$ smooth and the square Cartesian (i.e. $A_B=A\times_{\widetilde B}B$). It follows that $A$ and $A_B\to B$ are also smooth. Let $J_{0B}$ denote the ideal sheaf of $X_0\subset A_B$
and consider the natural commutative diagram
$$\xymatrix{J_0/J_0^2\ar[d]\ar[r]&J_{0B}/J_{0B}^2\ar[d]\\
\Omega_A|_{X_0}\ar[r]&**[r]\Omega_{A/\widetilde B}|_{X_0}\cong\Omega_{A_B/B}|_{X_0}.}$$
By definition, $\LL_{X_0/B}$ is the complex $J_{0B}/J_{0B}^2\to \Omega_{A_B/B}|_{X_0}$
and the diagram can be viewed as a morphism
\begin{equation}\label{eqn:LLrel}
\LL_{X_0}\to \LL_{X_0/B}.
\end{equation}
Since the natural map $J_0/J_0^2\to I$ factors through $J_{0B}/J_{0B}^2$, the Kodaira--Spencer class $\kappa(X_0/X)$ defined by (\ref{dg:KS}) factors naturally via
(\ref{eqn:LLrel}) giving a commutative diagram
\begin{equation} \label{bonn}
\xymatrix{\LL_{X_0}\ar[dr]_{\kappa(X_0/X)}
\ar[rr]&&\LL_{X_0/B}\ar[dl]^{\kappa(X_0/X/B)}\\
&I[1].\!\!&}
\end{equation}
We call $\kappa(X_0/X/B)\colon\LL_{X_0/B}\to I[1]$ the \emph{relative truncated Kodaira--Spencer class} of $X_0\subset X$.

\smallskip

Similarly, the (absolute) truncated universal Atiyah class $\alpha_{X_0}:\ko_{\Delta_{X_0}}\to
i_{\Delta_{X_0}*}(\LL_{X_0})[1]$ can be composed with the push-forward
$i_{\Delta_{X_0}*}(\LL_{X_0})\to i_{\Delta_{X_0}*}(\LL_{X_0/B})$ of (\ref{eqn:LLrel})
to yield the \emph{relative truncated universal Atiyah class}
$$\alpha_{X_0/B}:\ko_{\Delta_{X_0}}\to
i_{\Delta_{X_0}*}(\LL_{X_0/B})[1].$$

The diagram
$$\xymatrix{&\ko_{\Delta_{X_0}}\ar[dr]^{\alpha_{X_0/B}}\ar[dl]_{\alpha_{X_0}}&\\
i_{\Delta_{X_0}*}(\LL_{X_0})[1]\ar[dr]_{\!\!\!\!\!\!i_{\Delta_{X_0}*}\kappa(X_0/X)~}\ar[rr]&& i_{\Delta_{X_0}*}(\LL_{X_0/B})[1]\ar[dl]^{\quad i_{\Delta_{X_0}*}\kappa(X_0/X/B)}\\
&i_{\Delta_{X_0}*}I[1]&}$$
is commutative: the top half by definition and the bottom by \eqref{bonn}.
Therefore the universal obstruction class $\varpi(X_0/X)$ can also be computed as the product of the relative
Atiyah class and the relative Kodaira--Spencer class:
\begin{equation}\label{eqn:univobstrrel}
\varpi(X_0/X)=i_{\Delta_{X_0}*}(\kappa(X_0/X/B))\circ\alpha_{X_0/B}.
\end{equation}

\section{Obstructions} \label{sect:Lieblich}

\subsection{Deformations and extension classes} By \eqref{eqn:coh0-1} we have the following natural
exact triangles on $X_0$:
\begin{eqnarray}
\tau^{\le-1}H\xymatrix{\ar[r]&} H\xymatrix{\ar[r]&}\ko_{\Delta_{X_0}},\quad\text{and} \label{Hadj} \\
\tau^{\le-2}H\xymatrix{\ar[r]&}\tau^{\le-1}H\xymatrix{\ar[r]&}
i_{\Delta_{X_0}*}I[1]. \label{QHadj}
\end{eqnarray}
The universal obstruction class of Definition \ref{uoc} is the composition
of $\ko_{\Delta_{X_0}}
\to\tau^{\le-1}H[1]$ from the first triangle \eqref{Hadj} and $\tau^{\le-1}H[1]\to
i_{\Delta_{X_0}*}I[2]$ from the second \eqref{QHadj}:
\begin{equation} \label{uoc2}
\varpi(X_0/X)\colon\ko_{\Delta_{X_0}}\xymatrix{\ar[r]&}
i_{\Delta_{X_0}*}I[2].
\end{equation}
Thinking of these as Fourier--Mukai kernels,
act them on a fixed perfect complex $E_0$ on $X_0$. The result is the exact
triangles
\begin{eqnarray}
Q_{E_0}\xymatrix{\ar[r]&} i^*i_*E_0\xymatrix{\ar[r]&} E_0,\quad\text{and} \label{Eadj} \\
(\tau^{\le-2}H)(E_0)\xymatrix{\ar[r]&}
Q_{E_0}\xymatrix@=25pt{\ar[r]^{\pi_{E_0}}&} E_0\otimes I[1].
\label{Qadj}
\end{eqnarray}
The triangle \eqref{Eadj} is adjunction $i^*i_*E_0\to E_0$, whose
kernel we denote $Q_{E_0}$. By \eqref{uoc2}, the composition of
$E_0\to Q_{E_0}[1]$ from the first triangle \eqref{Eadj} and
$\pi^{}_{E_0}\colon Q_{E_0}[1]\to E_0\otimes I[2]$ from the second
\eqref{Qadj} gives the \emph{obstruction class} of $E_0$,
\begin{equation} \label{Eoc}
\varpi(X_0/X)(E_0)\colon E_0\xymatrix{\ar[r]&} E_0\otimes I[2].
\end{equation}
In this section we justify this terminology.
\medskip

By a deformation of $E_0$ over $X$ we mean a perfect complex $E$
on $X$ whose derived restriction $i^*E$ is isomorphic to $E_0$.
Such a deformation sits in an obvious exact triangle
\begin{equation} \label{def}
i_*(E_0\otimes I)\xymatrix{\ar[r]&} E\xymatrix{\ar[r]&}
i_*E_0\xymatrix{\ar[r]^e&} i_*(E_0\otimes I)[1],
\end{equation}
defining, and defined by, an extension class
$e\in\Ext^1_X(i_*E_0,i_*(E_0\otimes I))$. Applying the derived
functor $i^*$ gives the horizontal triangle in the diagram
\begin{equation} \label{Qdiag}
\xy {\ar (0,13)*+{i^*i_*(E_0\otimes I)}; (24,13)*+{i^*E}}; {\ar
(28,13); (37,13)}; {\ar (44,13)*+{i^*i_*E_0};
(72,13)*+{i^*i_*(E_0\otimes I)[1]}^(.35){i^*e}};
(44,26)*{Q_{E_0}}; (72,0)*{E_0\otimes I[1].};
(67,21)*{{}_{\Psi_e}}; (50,26); (79,4); **\crv{(102,20)}; {\ar
(44,23); (44,16)}; {\ar (44,10); (44,0)*+{E_0}}; {\ar (28,10);
(41,3)_r}; {\ar (50,10); (63,3)_{\tilde e}}; {\ar (48,23);
(63,16)^(.65){\Phi_e}}; {\ar (71,10); (71,3)}; {\ar (80.2,5);
(78,3.2)};
\endxy
\end{equation}
The two lower vertical arrows are adjunction. Therefore the
central vertical triangle is \eqref{Eadj}, and the map $\tilde
e\in\Ext^1_{X_0}(i^*i_*E_0,E_0\otimes
I)\cong\Ext^1_{X}(i_*E_0,i_*(E_0\otimes I))$ is what corresponds
to $e$ under adjunction.

Since the map denoted $r$ (the restriction map) is an isomorphism, $\Phi_e$
is too. Conversely, given any extension $e\colon i_*E_0\to i_*(E_0\otimes
I)[1]$, the map $r$ is an isomorphism if $\Phi_e$ is. Therefore

\begin{lem} \label{lemEE}
Deformations $E$ of $E_0$ are equivalent to extensions
$e\in\Ext^1_X(i_*E_0,i_*(E_0\otimes I))$ for which $\Phi_e$
\eqref{Qdiag} is an isomorphism. $\hfill\square$
\end{lem}

\begin{lem} \label{easyhalf}
If $e$ describes a deformation $E$ of $E_0$ then the map
$\Psi_e$ of \eqref{Qdiag} is
$$
\pi^{}_{E_0}\colon Q_{E_0}\xymatrix{\ar[r]&} E_0\otimes I[1]
$$
defined in \eqref{Qadj}.
\end{lem}

\begin{proof}
We begin by stating two results whose verification is left
to the reader.

(i) Fix a morphism $f:Z\to Y$ and  complexes of sheaves  $\kf,\,\kh$ on $Z$
and $\kg$ on $Y$. There are the following compatibilities between
the adjunction $f^*f_*\to{\rm id}$ and the projection formula.
$$\spreaddiagramrows{-.5pc}
\xymatrix{f^*f_*\kf\otimes f^*\kg\ar[d]^\cong\ar[rr]^(.55){{\rm adj}\
\otimes\ \id}\ar@{}[rrdd]|{\quad\circlearrowleft}
&&\kf\otimes f^*\kg\ar@{=}[dd]\\
f^*(f_*\kf\otimes\kg)\ar[d]^\cong&\\
(f^*f_*)(\kf\otimes f^*\kg)\ar[rr]^(.55){{\rm adj}}&&\kf\otimes f^*\kg}
\qquad \xymatrix{f_*(f^*f_*\kf\otimes\kh)\ar[d]^\cong\ar[rr]^(.53)
{f_*({\rm adj}\ \otimes\ \id)}\ar@{}[rrdd]|\circlearrowleft
&&f_*(\kf\otimes\kh)\ar@{=}[dd]\\
f_*\kf\otimes f_*\kh\ar[d]^\cong&\\
f_*(\kf\otimes f^*f_*\kh)\ar[rr]^(.53){f_*(\id\ \otimes\ {\rm
adj})}&&f_*(\kf\otimes\kh)}$$

(ii) Suppose $\kf$ is a Fourier--Mukai kernel on $X_0\times X_0$
for a functor $\Phi_\kf$ going from the second factor to the
first\footnote{In this proof only,
for notational reasons we write the kernels for Fourier--Mukai
functors from $Y$ to $Z$ as objects on $Z\times Y$; the reverse
of our usual convention. Therefore, for instance,
$h_*\ko_{\Delta_{X_0}}$ is the kernel for the functor $i^*$ (not
$i_*$, as it would be with the other convention).}
and $\kg=h_*\kf$ is its direct image under $h:X_0\times
X_0\ \hookrightarrow\,X_0\times X$ viewed as Fourier--Mukai kernel
for a functor $\Phi_\kg$ going from $X$ to $X_0$. Then
$\Phi_\kg\cong\Phi_\kf\circ i^*$.
\medskip

Start with the trivial case $E_0=\ko_{X_0}$ with canonical
deformation $E=\ko_X$. Then
it is easy to see that $Q_{\ko_{X_0}}\cong i^*i_*I[1]$ and both
maps $\Psi_{e}$ and $\pi^{}_{\ko_{X_0}}$ are the adjunction
\begin{equation} \label{Q0}
i^*i_*I[1]\xymatrix{\ar[r]&} I[1],
\end{equation}
i.e. the map given by taking $\kh^{-1}$.

For general $E_0$, we tensor the diagram \eqref{Qdiag} for the
deformation $\ko_X$ of $\ko_{X_0}$ by $E_0$, and use $E$ to pass
through $i_*$s and $i^*$s; for instance
$$
E_0\otimes i^*i_*I\ \cong\ i^*(E\otimes i_*I)\ \cong\ i^*i_*(E_0\otimes I).
$$
Using (i) we find the adjunction maps also survive tensoring with $E_0$,
so the result is the diagram \eqref{Qdiag} for the deformation $E$
of $E_0$.

So we are left with showing that $\pi^{}_{E_0}$ in \eqref{Qadj} is
$E_0$ tensored with \eqref{Q0} when $E_0=i^*E$. Using $E$ to
commute $E_0$ past $i^*i_*$ again, we find it is sufficient to
show that $\pi^{}_{E_0}$ is the adjunction $i^*i_*(E_0\otimes
I)[1]\to E_0\otimes I[1]$ when $E_0=i^*E$. \medskip

We do this at the level of Fourier--Mukai transforms.
Recall that the Fourier--Mukai kernel for $i^*i_*$ is
$H=h^*h_*\ko_{\Delta_{X_0}}$, so by (ii) the composition
$i^*i_*i^*$ has kernel
$h_*h^*h_*\ko_{\Delta_{X_0}}$ on $X_0\times X$. The second
diagram in (i) yields the commutativity of
\spreaddiagramrows{-.5pc}
$$\xymatrix{ h_*h^*h_*\ko_{\Delta_{X_0}} \dto \ar[r]^-\sim &
h_*((h^*h_*\ko_{X_0\times X_0})\otimes
\ko_{\Delta_{X_0}}) \dto \\
h_*\ko_{\Delta_{X_0}} \ar[r]^-\sim&  h_*(\ko_{X_0\times
X_0}\otimes \ko_{\Delta_{X_0}}),}
$$
with  the vertical arrows the natural adjunctions, which coincide
as maps to $h_*\ko_{\Delta_{X_0}}$ -- they are both the
$\kh^0$-map.
 Taking kernels of the vertical maps gives
$h_*(\tau^{\le-1}H)$ on the left hand side. On the right hand side
resolving $h_*\ko_{X_0\times X_0}$ by $h_*(\ko_{X_0}\boxtimes
I)\to\ko_{X_0\times X}$ shows that the kernel is isomorphic to
$$
h_*(h^*h_*(\ko_{X_0}\boxtimes I[1])\otimes\ko_{\Delta_{X_0}})=
h_*h^*h_*(i_{\Delta_{X_0}*}(I)[1]).
$$
Therefore $h_*(\tau^{\le-1}H)\cong h_*h^*h_*(i_{\Delta_{X_0}*}(I)[1])$. Applied
to any perfect complex $E$ on $X$ this gives the isomorphism
$Q_{E_0}\cong i^*i_*(E_0\otimes I)[1],$ where $E_0=i^*E$.

Now take the $\kh^{-1}$-map
$$
h_*h^*h_*(i_{\Delta_{X_0}*}(I)[1])\xymatrix{\ar[r]&}
h_*(i_{\Delta_{X_0}*}(I)[1]).
$$
Applied to $E$ this induces the adjunction $i^*i_*(E_0\otimes
I)[1]\to E_0\otimes I[1]$, where $E_0=i^*E$. But by construction it is the
pushforward by $h_*$ of the second arrow of \eqref{QHadj}, and so applied
to $E$ it gives $\pi^{}_{E_0}$ by its definition \eqref{Qadj}.
\end{proof}

This section is devoted to proving the converse. Together with
Lemma \ref{easyhalf} this proves

\begin{thm} \label{biggy}
Deformations $E$ of $E_0$ are equivalent to extensions
$e\in\Ext^1_X(i_*E_0,i_*(E_0\otimes I))$ for which
$\Psi_e=\pi^{}_{E_0}$.
\end{thm}

This implies the main Theorem of the Introduction:

\begin{cor} \label{bigtoo}
There is a deformation $E$ of $E_0$ if and only if
$$
0=(\id_{E_0}\otimes\kappa(X_0/X))\circ A(E_0)\ \in\
\Ext^2_{X_0}(E_0,E_0\otimes I),
$$
in which case the deformations form a torsor over $\Ext^1_{X_0}
(E_0,E_0\otimes I)$.
\end{cor}

\begin{proof}
Apply $\Hom_{X_0}(\ \cdot\ ,E_0\otimes I[1])$ to the exact triangle
$Q_{E_0}\to i^*i_*E_0\to E_0$:
\begin{equation} \label{EEE}
\spreaddiagramcolumns{-.8pc}
\spreaddiagramrows{-1pc}
\xymatrix{
& \Ext^1_X(i_*E_0,i_*(E_0\otimes I)) \ar@{=}[d] \\
\Ext^1_{X_0}(E_0,E_0\otimes I)\rto& \Ext^1_{X_0}(i^*i_*E_0,E_0\otimes I)\rto&
\Ext^1_{X_0}(Q_{E_0},E_0\otimes I) \rto& \Ext^2_{X_0}(E_0,E_0\otimes I).}
\end{equation}
Comparing to \eqref{Qdiag}, the second arrow takes an extension $\tilde e$
to the map $\Psi_e$. So Theorem \ref{biggy} says that $\tilde e$ corresponds
to a deformation $E$ of $E_0$ if and only if $\Psi_e=\pi^{}_{E_0}$.

The last arrow in \eqref{EEE} maps $\pi^{}_{E_0}$ to its
composition with $E_0\to Q_{E_0}[1]$, which by \eqref{Eoc} is our
obstruction class $\varpi(E_0)$. By Theorem \ref{product4} applied
to $E_0$ this is the product $(\id_{E_0}\otimes\kappa(X_0/X))\circ
A(E_0)$. Therefore the first assertion of the  Corollary follows
from the exact sequence \eqref{EEE}. For the second part, observe
that if a deformation $E$ exists then the adjunction $i^*i_*E_0\to
E_0$ is split (by the inverse of the map $r$ in \eqref{Qdiag}) and hence
$\Ext^1_{X_0}(E_0,E_0\otimes I)\to
\Ext^1_{X_0}(i^*i_*E_0,E_0\otimes I)$ is injective.
\end{proof}

By Lemmas \ref{lemEE} and \ref{easyhalf}, Theorem \ref{biggy} is
equivalent to the following assertion. For a given
$e\in\Ext^1_X(i_*E_0,i_*(E_0\otimes I))$,
 \smallskip
\begin{center}\emph{$\Phi_e$ is an isomorphism if and only if
$\Psi_e$  is the morphism $\pi^{}_{E_0}\colon Q_{E_0}\to
E_0\otimes I[1]$.}\end{center}

\smallskip\noindent Proving the `if' part of this directly would give the most
satisfactory treatment of the obstruction theory via the methods
of this paper. This seems to be extremely difficult, though we are
not entirely sure why.

\subsection{Proof in the affine case} \label{sec:affine}

Instead we mimic Lieblich, passing to the local case first.  So we
prove Theorem \ref{biggy} for $X_0\subset X$ both affine as
schemes (and not just as schemes over $B$), where we may assume
that $E_0$ is a bounded complex $E_0^\bullet$ of free modules.
Where possible we use similar notation to Lieblich in our slightly
different setting.

So we fix an extension $e$ with $\Psi_e=\pi^{}_{E_0}$. We observed in
the proof of Corollary \ref{bigtoo} that the image of $\pi_{E_0}^{}$ in the
exact sequence \eqref{EEE} is $\varpi(E_0)$, so $\varpi(E_0)=0$ by assumption.

Thinking of the complex $I_{\Delta_{X_0}\subset X_0\times X}\to\ko_{X_0\times
X}$ as a Fourier--Mukai kernel, quasi-isomorphic to the kernel $\ko_{\Delta_{X_0}}$
of the pushforward $i_*$, we apply it to $E_0^\bullet$. This gives a
resolution
\begin{equation} \label{bullet}
\xymatrix{
\Big(K^\bullet\rto & \Gamma^\bullet\Big) \rto^{\cong} & i_*E_0^\bullet}
\end{equation}
which moreover is a quasi-isomorphism for each fixed index $\bullet$. Here
$\Gamma^j:=
\Gamma(E_0^j)\otimes_{\ko(B)}\ko_X$ and $K^j$ is the kernel of the evaluation
map $\Gamma^j\to i_*E_0^j$.

Pulling the Fourier--Mukai kernels back to $X_0\times X_0$ gives
the exact sequence \eqref{t-H} with extension class
$\varpi(X_0/X)$ in \eqref{uoc2}. Applying this to $E_0^\bullet$
shows that pulling \eqref{bullet} back to $X_0$ gives the
compatible long exact sequences
\begin{equation} \label{LES}
\xymatrix{
0 \rto& E_0^\bullet\otimes I \rto& \K^\bullet|_{X_0} \rto& \Gamma^\bullet|_{X_0}
\rto& E_0^\bullet \rto& 0,}
\end{equation}
with extension class the obstruction class $\varpi(E_0)\in
\Ext^2_{X_0}(E_0,E_0\otimes I)$. (We use $\ \cdot\ |_{X_0}$ to denote the
the \emph{un}derived functor $i^*$.)

Pick a free deformation $E^j$ of each (free) $E_0^j$ to $X$. Since $\Gamma^j$
is free, we can split it as $N^j\oplus E^j$ so that
the map $\Gamma^j\to i_*E_0^j$ of \eqref{bullet} is projection to the second
factor of $N^j\oplus E^j$ followed by restriction to $X_0$. Being in the
kernel of this map, $N^j$ lifts naturally to $K^j$ in \eqref{bullet}. The
result is the following
non-canonical splitting of \eqref{bullet}, incompatible with the differential.
\begin{equation} \label{noncanon}
\xymatrix{
0 \rto& N^j\oplus i_*(E_0^j\otimes I) \rto& N^j\oplus E^j \rto& i_*E_0^j
\rto& 0.}
\end{equation}
Restricting to $X_0$ splits \eqref{LES} as
\begin{equation} \label{mess}
\spreaddiagramrows{-.7pc}
\spreaddiagramcolumns{-.4pc}
\xymatrix{
0 \rto& E_0^j\otimes I \rto\dto& N^j|_{X_0}\oplus(E_0^j\otimes I) \rto\dto&
N^j|_{X_0}\oplus E^j_0 \rto\dto& E_0^j \rto\dto& 0 \\
0\rto&E_0^{j+1}\otimes I\rto\dto&N^{j+1}|_{X_0}\oplus(E_0^{j+1}\otimes I)
\rto\dto&N^{j+1}|_{X_0}\oplus E^{j+1}_0\rto\dto&E_0^{j+1}\rto\dto&0 \\
0\rto&E_0^{j+2}\otimes I \rto& N^{j+2}|_{X_0}\oplus(E_0^{j+2}\otimes I)\rto&
N^{j+2}|_{X_0}\oplus E^{j+2}_0 \rto& E_0^{j+2} \rto& 0.\!}
\end{equation}

With respect to this splitting the differential $\Gamma^j|_{X_0}\to\Gamma^{j+1}|_{X_0}$
takes the form
$$
\left(\begin{array}{cc} \!\!* & \beta \\ \!\!* & d_{E_0}\!\!\end{array}\right)
\!\!\! \begin{array}{c} \\ . \end{array}
$$
Since $N^{j+1}|_{X_0}$ lies naturally in $K^{j+1}|_{X_0}$, $\beta$ defines
a map from $E_0^j$ to $K^{j+1}|_{X_0}$. Composing with the differential to
$K^{j+2}|_{X_0}$ and projecting to $E_0^{j+2}\otimes I$ defines a map $E_0\to
E_0\otimes I[2]$ that represents the extension class $\varpi(E_0)$ of \eqref{LES}.

We now describe this map on $X$ using \eqref{noncanon} in place of \eqref{mess}
on $X_0$. The differential on $\Gamma^j$ splits as
\begin{equation} \label{matrix}
\left(\begin{array}{cc} \!\!* & \gamma \\ \!\!* & d_E\!\!\end{array}\right)
\!\!\! \begin{array}{c} \\ , \end{array}
\end{equation}
where $\gamma$ factors through a map from $E^j$ to $K^{j+1}$. (After restricting
to $K^{j+1}|_{X_0}$ this vanishes on $E^j\otimes i_*I$
and descends to the map $\beta\colon E_0^j\to K^{j+1}|_{X_0}$ above.)
The splitting \eqref{matrix}
defines $d_E$ on $E^\bullet$, which need not square to zero, though
it covers the genuine differential $d_{E_0}$.

The inclusion $K^\bullet\ \mono\,\Gamma^\bullet$ commutes with the
differentials $d_K,d_\Gamma$, so instead of applying $d_K$ as
above we apply $d_\Gamma$ to $\gamma$. Applying the projection $P_E$ to the
$E^{j+2}$ component of $\Gamma^{j+2}$
gives a map $P_E(d_\Gamma\gamma)\colon i_*E^j_0\to E^{j+2}$. This
projects to zero in $i_* E_0^{j+2}$, so has image in
$E^{j+2}\otimes i_*I$. By its construction it is $i_*$ applied to
the map $\varpi(E_0)\colon E_0^j\to E_0^{j+2}\otimes I$ described
above.

However, the differential on $\Gamma^\bullet$ squares to zero, so applying
this to $E^j\subset\Gamma^j$ gives, by \eqref{matrix},
$$
d_E^2+P_E(d_\Gamma\gamma)=0.
$$
Thus our obstruction $\varpi(E_0)$ is precisely $[-d_E^2]$. By
assumption this is zero in $\Ext^2_{X_0}(E_0,E_0\otimes I[2])$, so
$-d_E^2=d_{E_0}h+hd_{E_0}$ for some homotopy $h\colon E^j_0\to
E^{j+1}_0\otimes I$. Lift $h$ to $\tilde h\colon E^j\to
E^{j+1}_0\otimes I$ and add to $d_E$ using the inclusion
$E^{j+1}_0\otimes I\ \mono\,E^{j+1}$. The resulting $d_E$ has
square zero and so defines a deformation $(E^\bullet,d_E)$ of
$E^\bullet_0$, corresponding to an extension
$e'\in\Ext^1_X(i_*E_0,i_*(E_0\otimes I))$.

Similarly any other differential on $E^\bullet$ (covering $d_{E_0}$
on $E_0^\bullet$) differs from $d_E$ by a map $f\colon E_0^j\to E_0^{j+1}\otimes
I$ satisfying $d_{E_0}\circ f+f\circ d_{E_0}=0$. Therefore $[f]\in
\Ext^1_{X_0}(E_0,E_0\otimes I)$, with the resulting
map from $\Ext^1_{X_0}(E_0,E_0\otimes I)$ to $\Ext^1_X(i_*E_0,i_*(E_0\otimes
I))$ being the first map in the exact sequence \eqref{EEE}. This describes
the fibre over $\pi^{}_{E_0}$ in \eqref{EEE}, so $e-e'$ is in the image of
some $[f]\in\Ext^1_{X_0}(E_0,E_0\otimes I)$.
Altering $d_E$ by $f$ as above gives the deformation
$(E^\bullet,d_E)$ of $E^\bullet_0$ corresponding to $e$, as required.

\subsection{Proof in the general case}
We can now prove Theorem \ref{biggy} for any $X_0\subset X$. Suppose that
$\varpi(E_0)=0$: i.e. the exact sequence of Fourier--Mukai kernels
\eqref{t-H} applied to $E_0$ gives an exact triangle with zero extension
class. Then we claim that the same is true on restriction to any
affine open subset $U$: the restriction of the exact triangle is the restriction
of \eqref{t-H} to $U\times U$ applied to $E_0|_U$. (This is not true for
general Fourier--Mukai kernels, but is true here since $\tau^{\ge-1}H$ is
quasi-isomorphic to a complex supported, set-theoretically, on the diagonal,
making it a local operator. Using this second complex as a kernel makes the
result clear.)

Therefore any $e\in\Ext^1_X(i_*E_0,i_*(E_0\otimes I))$ mapping to
$\pi^{}_{E_0}$ in \eqref{EEE} defines an extension $E$ which, on
restriction to $U$, is a deformation of $E_0$ by the previous
Section \ref{sec:affine}. Therefore it is a deformation of $E_0$
on $X$ itself, since the condition that the map $r\colon i^*E\to
E_0$ of \eqref{Qdiag} be a quasi-isomorphism is local.

\subsection{Relative version.} We will need the following  relative version of Corollary \ref{bigtoo}. First,
in the situation of Section \ref{sec:relKS}, one defines the \emph{relative truncated Atiyah class} as
\begin{equation} \label{relAt}
A(E_0/B):=\alpha_{X_0/B}(E_0)\in\Ext^1_{X_0}(E_0,E_0\otimes \LL_{X_0/B})
\end{equation}
or, equivalently, as the image of $A(E_0)$ under
the projection $$\Ext^1_{X_0}(E_0,E_0\otimes \LL_{X_0})\to \Ext^1_{X_0}(E_0,E_0\otimes \LL_{X_0/B}).$$
Applying the identity $\varpi(X_0/X)=i_{\Delta_{X_0}*}(\kappa(X_0/X))\circ\alpha_{X_0}=i_{\Delta_{X_0}*}(\kappa(X_0/X/B))\circ
\alpha_{X_0/B}$ of \eqref{eqn:univobstrrel} to $E_0$ gives
\begin{equation} \label{relobs}
\varpi(X_0/X)(E_0)=(\id_{E_0}\otimes\kappa(X_0/X))\circ A(E_0)=
(\id_{E_0}\otimes\kappa(X_0/X/B))\circ A(E_0/B).
\end{equation}
Thus, there is a deformation $E$ of $E_0$ if and only if
$$0=(\id_{E_0}\otimes\kappa(X_0/X/B))\circ A(E_0/B)\ \in \Ext^2_{X_0}(E_0,E_0\otimes I).
$$


\section{Virtual cycles} \label{sect:appl}

\subsection{Setup}
In this Section we work over a base $B$ defined
over a field $k$ whose characteristic does not
divide the rank of the complexes considered below.

Fix an $n$-dimensional smooth projective connected morphism $X\to
B$ and a line bundle $\kl$ over $X$. Denote by $i_b\colon X_b\
\mono\,X$ the fibre of $X$ over a closed point $b\in B$. Given a
separated algebraic space $\km/B$ and a point closed point
$m\in\km$ sat over $b$, we denote by $i_m\colon X_b\
\mono\,X\times_B\km$ the inclusion of $X_b\times\{m\}$.

Suppose that $\km/B$ is a \emph{relative fine moduli space of perfect simple
complexes} with fixed determinant $\kl$ on $X/B$. By this we mean that there
is a perfect complex
\begin{equation} \label{univ}
\EE\in D^b(X\times_B\km)
\end{equation}
(the \emph{universal complex\footnote{This can also be done with a twisted
universal complex, though we do not need such generality for the application
\cite{PT} we have in mind.}} for the moduli space
) with $\det(\EE)\cong\pi_\km^*M\otimes\pi_X^*\kl$ for some
$M\in\Pic(\km)$ such that $\km/B$ parameterises the complexes
$i_m^*\EE$ (simple, of determinant $i_b^*\kl$) on the fibres
$X_b\times\{m\}$ in the following sense:
\begin{itemize}
\item the set of morphisms of $B$-schemes $f\colon S\to\km$, and
\item the set of equivalence classes of perfect complexes $E$ over
$X\times_BS$ whose restriction $i_s^*E$ to any fibre $X_s,\ s\in
S$, is isomorphic to $i_m^*\EE$ for some $m\in\km$ and such that
$\det(E)\cong\pi_S^*M_S\otimes \pi_X^*\kl$ for some
$M_S\in\Pic(S)$
\end{itemize}
are put in bijection by assigning to $f$ the perfect complex
$(\id\times f)^*\EE$ over $X\times_B S$. By definition, two
complexes on $X\times_BS$ are equivalent if they differ by a twist
with a line bundle coming from $S$. (This definition ensures that
$\km$ is locally complete, in the sense that any deformation of a
complex in $\km$ is also in $\km$.)

As Lieblich rightly pointed out, it is most natural to work with Artin
stacks at this point. In \cite{Lieb} it was shown under very
general assumptions that complexes as objects in the derived category
form an Artin stack which is locally of finite type. However, in
the applications we have in mind, the moduli space is also known to exist
in the above sense as a quasi-projective variety or an
algebraic space. In order to keep the discussion simple (and due
to our incompetence), we decided to keep to this setting.

\medskip

In the following we shall tacitly assume that $\rk(\EE)$ is a
constant (not divisible by the characteristic of our base field
$k$ as mentioned above). Also, we will assume that the Euler
characteristic $\chi(i_m^*\EE,i_m^*\EE)$ is independent of $m\in
\km$, for instance by fixing the numerical invariants of the
complexes we are interested in. Finally we will assume that $\km$ admits
an embedding into a smooth (over $B$) ambient space $A$, though in Section
\ref{remove} below we describe how to remove this assumption.


\subsection{Relative obstruction theory} \label{rel}

Let $\pi_X$ and $\pi_\km$ denote the projections from $X\times_B\km$
to $X$ and $\km$ respectively. For simplicity we will also assume
that $\rk(\EE)\ne0$ for now; we deal with the rank zero case at
the end of this section. Then since the composition of $\id\colon
\ko_{X\times_B\km}\to\sHom(\EE,\EE)$ and $\tr\colon\sHom(\EE,\EE)
\to\ko_{X\times_B\km}$ is multiplication by $\rk(\EE)$ we get a
splitting
$$
\sHom(\EE,\EE)\ \cong\ \sHom(\EE,\EE)_0\oplus\ko_{X\times_B\km}.
$$

Assuming $\km$ admits a smooth embedding, we have defined the truncated Atiyah
class \eqref{eq:truncAt} of the universal complex $\EE$,
$$
A(\EE)\in\Ext^1_{X\times_B\km}(\EE,\EE\otimes\LL_{X\times_B\km}).
$$
Via $\LL_{X\times_B\km}\to\LL_{X\times_B\km/X}=\pi_\km^*\LL_{\km/B}$ this maps to the relative truncated Atiyah class \eqref{relAt}
$$
A(\EE/X)\in\Ext^1_{X\times_B\km}(\sHom(\EE,\EE)_0,\pi_\km^*\LL_{\km/B}).
$$
By Verdier duality along the projective morphism $\pi_\km$ this is isomorphic
to
$$
\Ext^{1-n}_{\km}(\pi_{\km*}(\sHom(\EE,\EE)_0\otimes\omega_{\pi^{}_\km}),\LL_{\km/B}),
$$
where $\omega_{\pi^{}_\km}=\pi_X^*\omega^{}_{X/B}$ is the relative dualizing
sheaf. We obtain a map
\begin{equation} \label{obsthy}
\pi_{\km*}(\sHom(\EE,\EE)_0\otimes\pi_X^*\omega^{}_{X/B})[n-1]\xymatrix{\ar[r]&}\LL_{\km/B}.
\end{equation}
The results of this paper have essentially proved the following.

\begin{thm} \label{BF}
In the notation of \cite{BF} (Definition 4.4 and Section 7),
the map \eqref{obsthy} is a relative obstruction theory for $\km$.
\end{thm}

\begin{proof}
Fix a morphism of $B$-schemes,
$$
f\colon S_0\xymatrix{\ar[r]&}\km,
$$
and an extension $S_0\subset S$ with ideal $I$ such that $I^2=0$.

The composition of the pullback
of \eqref{obsthy},
\begin{equation} \label{fobs}
f^*\pi_{\km*}(\sHom(\EE,\EE)_0\otimes\pi_X^*\omega^{}_{X/B})[n-1]\xymatrix{\ar[r]&}
f^*\LL_{\km/B}
\end{equation}
with the natural map $f^*\LL_{\km/B}\to\LL_{S_0/B}$ followed by the relative
truncated Kodaira--Spencer class \eqref{bonn} of $S_0\subset S$,
$$
\kappa(S_0/S/B)\in\Ext^1_{S_0}(\LL_{S_0/B},I),
$$
gives an element
\begin{equation} \label{oext2}
o\in\Ext^{2-n}_{S_0}(f^*\pi_{\km*}(\sHom(\EE,\EE)_0\otimes\pi_X^*\omega^{}_{X/B}),I).
\end{equation}
By \cite[Thm.\ 4.5]{BF} we must show that $o$ vanishes if and only
if there exists an extension from $S_0$ to $S$
of the map $f$ of
$B$-schemes, and that when $o=0$ the set of extensions is a torsor
under $\Ext^{1-n}(f^*\pi_{\km*}(\sHom(\EE,\EE)_0\otimes
\pi_X^*\omega^{}_{X/B}),I)$. (In fact \cite[Thm.\ 4.5]{BF} only proves
this in the absolute, not relative, case. The proof works just the
same over a base.)

Let
$$
\bar f=\id\times f\colon X\times_B S_0\xymatrix{\ar[r]&} X\times_B\km
$$
and let $\bar \pi$ be the projection
$$
\bar{\pi}\colon X\times_B S_0\xymatrix{\ar[r]&} S_0.
$$
Since $\bar\pi$ is flat, the composition of \eqref{fobs} with
$f^*\LL_{\km/B}\to\LL_{S_0/B}$ is
$$
\pi_{\km*}(\sHom(\bar f^*\EE,\bar
f^*\EE)_0\otimes\pi_X^*\omega^{}_{X/B})[n-1]\xymatrix{\ar[r]&}\LL_{S_0/B}.
$$
By the functoriality of truncated Atiyah classes this is the trace-free part of the 
relative truncated Atiyah class $A(\bar f^*\EE/X)$. Therefore its
product with the relative truncated Kodaira--Spencer class
$$\kappa(X\times_BS_0/X\times_BS/X)=\bar\pi^*\kappa(S_0/S/B),$$
is by \eqref{relobs} precisely the trace-free part of our obstruction class $\varpi(X\times_BS_0/X\times_B S)(\bar f^*\EE)$.
We denote this by
\begin{equation} \label{last}
o\in\Ext^{2-n}_{S_0}(\bar\pi_*(\sHom(\bar f^*\EE,\bar
f^*\EE)_0\otimes\omega_{\bar\pi}),I) \cong\Ext^2_{X\times_B
S_0}(\bar f^*\EE,\bar f^*\EE\otimes\bar \pi^*I)_0,
\end{equation} the isomorphism being Verdier duality for $\bar\pi$.

That the trace of the obstruction class is the obstruction to
deforming $\det(\bar f^*\EE)$ is classical. For instance, this is
proved for locally free sheaves in \cite[Prop.\ 3.15]{Th}, and
then, using this, for complexes of locally free sheaves in
\cite[Thm.\ 3.23]{Th}. Since $\det(\bar f^*\EE)$ deforms as $\kl$,
the trace of $\varpi(E_0)$ vanishes. Similarly the choices of
extension are governed by the trace-free part of $\Ext^1$
\cite[Thm.\ 3.23]{Th}. A better, more thorough account using the language of this paper is now available in \cite{La}.

Thus, by Corollary \ref{bigtoo}, $o=0$ if and only if there exists
a deformation of $\bar f^*\EE$ from $X\times_BS_0$ to
$X\times_BS$, in which case the deformations with fixed determinant
form a torsor under $\Ext^1_{X\times_BS_0}(\bar f^*\EE,\bar
f^*\EE\otimes\bar\pi^*I)_0$.

But deformations of $\bar f^*\EE$ from $X\times_BS_0$ to $X\times_BS$ are in one-to-one correspondence with extensions from $S_0$ to $S$ of
the $B$-map $f$, by the definition of a relative fine moduli space \eqref{univ}.
So we are done.
\end{proof}

\subsection{Virtual cycles} \label{virtcy}
To get a virtual cycle we must make this obstruction theory
\emph{perfect}, or 2-term. The simplest way to ensure this is to
demand that each simple complex $E=i_m^*\EE$ in $\km$ satisfies
\begin{equation} \label{conditions}
\Ext^i_{X_b}(E,E)_0=0,\quad i\ne1,2.
\end{equation}
The following is then a standard consequence of the
Nakayama Lemma, but we give the argument in full.

\begin{lem} \label{perf}
Under the conditions \eqref{conditions}, the complex
$\pi_{\km*}(\sHom(\EE,\EE)_0)$ on $\km$ is quasi-iso\-morphic to a 2-term
complex of locally free sheaves $(F^1\to F^2)$ in degrees 1 and 2.
\end{lem}

\begin{proof}
Resolve $\sHom(\EE,\EE)_0\to A^\bullet$ by a finite complex of locally free
sheaves $A^j$ with vanishing $R^i\pi_{\km*}(A^j)$ for $i\ne0$ and all $j$.
This can
be achieved by picking a very negative locally free resolution of the derived
dual, dualising this and then truncating at some $j\gg0$.

Thus for each $j$, $\pi_{\km*}(A^j)$ is canonically quasi-isomorphic
to some locally free sheaf $F^j$. The resulting complex
$F^\bullet$ is a representative of $\pi_{\km*}(\sHom(\EE,\EE)_0)$.

By base change, the restriction of $F^\bullet$ over a closed point
$m\in\km$ (over $b\in B$) is a complex of vector spaces computing
$\Ext^*(i_m^*\EE,i_m^*\EE)_0$. Therefore by \eqref{conditions},
the restriction of $F^\bullet$ to each fibre $X\times\{m\}$ has
cohomology only in degrees 1 and 2. So if $F^{j>2}$ is the last
nonzero term of $F^\bullet$, then on each fibre
\begin{equation}\label{xdee}
F^{j-1}\xymatrix{\ar[r]&} F^j
\end{equation}
is surjective. Hence, the map \eqref{xdee} is surjective globally
with locally free kernel. Replacing $F^{j-1}$ by the kernel
and $F^j$ by zero, we can inductively assume that $j=2$.
Similarly if the first nonzero term is $F^{i<1}$,
then $F^i\to F^{i+1}$ is injective on fibres with
locally free cokernel.  We conclude that $F^\bullet$ is quasi-isomorphic
to a 2-term complex $(F^1\to F^2)$.
\end{proof}

\begin{cor}
The obstruction theory \eqref{obsthy} is \emph{perfect} in the
sense of \cite[Def.\ 5.1]{BF}. Each $\km_b$, $b\in B$ carries a
virtual cycle $[\km_b]^{vir}$ of dimension $\chi(\kh
om(i_m^*\EE,i_m^*\EE)_0)$ ($m\in\km$) that is deformation
invariant: there exists a cycle $[\km]^{vir}$ on $\km$ with
$i_b^![\km]^{vir}=[\km_b]^{vir}$.
\end{cor}

\begin{proof}
By Verdier duality and Lemma \ref{perf}, the complex
$\pi_{\km*}(\sHom(\EE,\EE)_0\otimes\pi_X^*\omega^{}_X)[n-1]$ of
\eqref{obsthy} is quasi-isomorphic to a 2-term complex of locally
free sheaves $(F^2)^\vee\to(F^1)^\vee$ in degrees $-2$ and $-1$.
The statements about virtual cycles then follow from \cite[Prop.\
7.2]{BF}.
\end{proof}

By Serre duality the conditions \eqref{conditions} cannot reasonably be satisfied
in dimensions $\ge4$, and are automatically satisfied in dimensions $\le2$
by simple complexes $E$ that have no negative Exts. So the critical dimension
is 3, where they can be replaced by the conditions that $E$ be simple and
$$
\begin{array}{ll}
\mathrm{(i)} & \Ext^i_{X_b}(E,E)_0=0,\quad i<0,\text{ and} \\
\mathrm{(ii)} &
H^0(X_b,\omega^{}_{X_b})\xymatrix{\ar[r]&}\Hom(E,E\otimes\omega^{}_{X_b})
\text{ is an isomorphism}.
\end{array}
$$
Serre duality then implies that $\tr\colon\Ext^3_{X_b}(E,E)\to
H^3(\ko_{X_b})$ is an isomorphism and hence $\Ext^i_{X_b}(E,E)_0$
vanishes for $i\ge3$.

For \emph{Calabi--Yau} threefolds
($\omega^{}_{X_b}\cong\ko_{X_b}$) we only need (i) because $E$
being simple implies (ii). In this case $\chi(\kh
om(i_m^*\EE,i_m^*\EE)_0)=0$ so for moduli of simple complexes
satisfying (i) we get a degree 0 virtual cycle in each $\km_b$.
Thus if $\km/B$ is proper then we get an invariant -- the length
of the virtual cycle -- which is independent of $b$. By duality
the obstruction theory \eqref{obsthy} is self-dual in the sense of
Behrend \cite{Behrend} so if $\km$ is quasi-projective then one can use his
weighted Euler
characteristic approach to the virtual class, and extend the
definition of invariants to non-proper $\km$ (at the expense of
losing deformation invariance).

Another situation where (ii) is automatic is for \emph{arbitrary}
threefolds and complexes satisfying
$$
\mathrm{(ii')} \quad \sHom^0(E,E)=\ko_{X_b} \quad\text{and}\quad
\cal{E}xt^i(E,E)=0,\quad
i<0. $$
This implies both (ii) and that $E$ is simple, and applies to the ideal sheaves
of \cite{MNOP, Th} and the complexes $I^\bullet$ of \cite{PT}.

\subsection{The rank zero case} \label{r0}
When $\rk(E)=0$ we consider simple complexes satisfying (ii) or
(ii$'$), but we do \emph{not} fix their determinant. So it is
natural to use
\begin{equation} \label{obsthy2}
\pi_{\km*}(\sHom(\EE,\EE)\otimes\pi_X^*\omega^{}_{X/B})[n-1]\xymatrix{\ar[r]&}\LL_{\km/B}.
\end{equation}
(cf. \eqref{obsthy}) as an obstruction theory. The left hand side is Verdier
dual to $\pi_{\km*}\sHom(\EE,\EE)[1]$, which by the methods of Lemma
\eqref{perf} is quasi-isomorphic to a complex of locally free sheaves
concentrated in degrees $-1$ to 2. Therefore this obstruction theory is not
\emph{perfect}.

To trim it down we first form
\begin{equation} \label{identity}
\mathrm{Cone}\,\big(\ko_{\km}\xymatrix{\ar[r]&}\pi_{\km*}\sHom(\EE,\EE)\big),
\end{equation}
where the arrow is the identity map. This map induces an
isomorphism $k\to\Hom^0(i_m^*\EE,i_m^*\EE)$ on every closed fibre
$X\times\{m\}$ of $\pi_\km$, so, by the same use of the Nakayama
Lemma as in the proof of Lemma \ref{perf}, the complex
\eqref{identity} is $\tau^{\ge1}\pi_{\km*}\sHom(\EE,\EE)$, a
complex of locally free shaves in degrees 1 to 3.

Similarly the top (degree 3) part of the trace map on
$\pi_{\km*}\sHom(\EE,\EE)$ lifts naturally to the $\kh^3$-map
$$
\tau^{\ge1}\pi_{\km*}\sHom(\EE,\EE)\xymatrix{\ar[r]&}
R^3\pi_{\km*}(\ko_{X\times_B\km})[-3],
$$
with kernel $\tau^{[1,2]}\pi_{\km*}\sHom(\EE,\EE)$. Now $R^3\pi_{\km*}
(\ko_{X\times_B\km})$ is locally free over $\km$, and after base-change to
the fibre over $m\in\km$ the trace map is an isomorphism $Ext^3(\EE_m,\EE_m)\to
H^3(\ko_{X_b})$ (Serre dual to the isomorphism (ii)). Therefore the standard
Nakayama Lemma
methods again show that $\tau^{[1,2]}\pi_{\km*}\sHom(\EE,\EE)$ is
quasi-isomorphic to a 2-term complex of locally free sheaves.

Finally the map \eqref{obsthy2} gives by Verdier duality the diagram
\begin{equation} \label{3mod}
\xymatrix{
(\pi_{\km*}\sHom(\EE,\EE))^\vee[-1] \rto & \LL_{\km/B} \ar@{<--}^a[ddl]+UR \\
(\tau^{\ge1}\pi_{\km*}\sHom(\EE,\EE))^\vee[-1] \uto\dto\urto \\
(\tau^{[1,2]}\pi_{\km*}\sHom(\EE,\EE))^\vee[-1].\!\!}
\end{equation}
The arrow labelled $a$ can be filled in uniquely since the cone on the lower
vertical arrow is $(R^3\pi_{\km*}(\ko_{X\times_B\km})[-3])^\vee$, a vector
bundle concentrated in degree $-3$, whereas $\LL_{\km/B}$ is concentrated in
degrees $-1$ and $0$.

Then we modify the proof of Theorem \ref{BF} to show that the map
$a$ is again a relative obstruction theory for $\km$ in the sense
of \cite[Def.\ 4.4]{BF}. To do this we must assume that the
$B$-schemes $S_0$ and $S$ in the proof of Theorem \ref{BF} are
\emph{affine} -- by the proof of \cite[Thm.\ 4.5]{BF} this is
sufficient to prove we get an obstruction theory. Following the
proof of Theorem \ref{BF} we get an obstruction $o$ not now in
$$
\HH^2(\bar\pi_*\sHom(\bar f^*\EE,\bar f^*\EE\otimes\bar
\pi^*I)_0)=\Ext^2(\bar f^*\EE,\bar f^*\EE\otimes\bar\pi^*I)_0
$$
as in (\ref{oext2}, \ref{last}) (where $\HH$ denotes hypercohomology),
but instead in
$$
\HH^2(\tau^{[1,2]}\bar\pi_*\sHom(\bar f^*\EE,\bar
f^*\EE\otimes\bar \pi^*I)).
$$
By the collapse of the Leray spectral sequence for $\bar\pi$ over affine
$S_0$, this is
\begin{eqnarray*}
H^0\big(\kh^2(\tau^{[1,2]}\bar\pi_*\sHom(\bar f^*\EE,\bar f^*\EE\otimes\bar
\pi^*I))\big) &=& H^0\big(\kh^2(\bar\pi_*\sHom(\bar f^*\EE,\bar f^*\EE\otimes\bar
\pi^*I))\big) \\ &=& \Ext^2(\bar f^*\EE,\bar f^*\EE\otimes\bar\pi^*I).
\end{eqnarray*}
In other words the modification \eqref{3mod} changed the obstruction theory
in degrees 0 and 3 only, so does not affect $\cal{E}xt^2$. Therefore the
same proof goes through as before and $a$
\eqref{3mod} yields our obstruction $\varpi(E_0)\in\Ext^2(\bar f^*\EE,\bar
f^*\EE\otimes\bar \pi^*I)$. In the same way the choices in $\Ext^1$ also
work out just as before.

\subsection{Removing the smooth embedding hypothesis} \label{remove}

For the applications in \cite{PT} $\km$ is projective so smoothly
embeddable. We also expect more general proper moduli spaces of
objects of the derived category of a smooth threefold to be
projective. However it is certainly easier to produce such moduli
spaces as algebraic spaces by abstract arguments starting from
Artin stacks. Therefore it would be nice to produce virtual cycles
without the assumption that $\km$ be smoothly embeddable.

This is certainly possible. Firstly we must produce the obstruction theory
\eqref{obsthy} on a general $\km$. One way to do this is to use Illusie's
full cotangent complex $L^\bullet_{\km/B}$
and Atiyah class. The working of Section \ref{rel} then proceeds exactly
as before without the embedding assumption, yielding
$$
\pi_{\km*}(\sHom(\EE,\EE)_0\otimes\pi_X^*\omega^{}_{X/B})[n-1]\xymatrix{\ar[r]&}L^\bullet_{\km/B}.
$$
Composing with the truncation $\tau^{\ge-1}\colon L^\bullet_{\km/B}\to\LL_{\km/B}$
recovers \eqref{obsthy}. Restricted to any subset of $\km$ which is smoothly
embeddable this gives the same map as before.

As used in Section \ref{r0} above, to check that
\eqref{obsthy} is an obstruction theory for $\km$
it is only necessary to consider \emph{affine} $B$-schemes $S_0$ and $S$
in the proof of Theorem \ref{BF}. Doing so, we only use their image in $\km$,
which is itself affine and therefore smoothly embeddable. Our theory therefore
applies on that image to prove Theorem \ref{BF}. The rest of the Section
is unaffected, so we obtain deformation invariant virtual cycles under the
same conditions as in Section \ref{virtcy}.


\end{document}